\documentclass{amsart}%
\usepackage{amsmath}
\usepackage{amsfonts}
\usepackage{amssymb}
\usepackage{graphicx}%
\providecommand{\U}[1]{\protect\rule{.1in}{.1in}}

\theoremstyle{definition}

\theoremstyle{remark}

\numberwithin{equation}{section}
\begin{document}
\title{Less Mundane Applications of the Most Mundane Functions}
\author{Pisheng Ding}
\address{Department of Mathematics, Illinois State University, Normal, Illinois 61790}
\email{pding@ilstu.edu}
\subjclass[2020]{26D99, 26B25, 52A40.}
\keywords{gradient, Cauchy-Schwarz inequality, GM-AM inequality, convex functions,
distance formula.}

\begin{abstract}
Linear functions are arguably the most mundane among all functions. However,
the basic fact that a multi-variable linear function has a constant gradient
field can provide simple geometric insights into several familiar results such
as the Cauchy-Schwarz inequality, the GM-AM inequality, and some distance
formulae, as we shall show.

\end{abstract}
\maketitle

\section{Preliminaries}

\noindent We set the stage for the applications by recalling some basics. For
$f(x,y)=ax+by$, $\nabla f\equiv a\mathbf{i}+b\mathbf{j}$. As a result, at any
point in $%
\mathbb{R}
^{2}$, it is in the direction $a\mathbf{i}+b\mathbf{j}$ that the directional
derivative attains its maximum value of $\left\vert \nabla f\right\vert
=\sqrt{a^{2}+b^{2}}$. Geometrically, this means that the slope of the planar
graph of $f$ equals $\left\vert \nabla f\right\vert =\sqrt{a^{2}+b^{2}}$. The
level sets of $f$ are the parallel lines perpendicular to $a\mathbf{i}%
+b\mathbf{j}$. Denote the level-$0$ line of $f$ by $L_{0}$. Any point in $%
\mathbb{R}
^{2}$ can be reached by traveling from a point on $L_{0}$ in either the normal
direction $a\mathbf{i}+b\mathbf{j}$ or its opposite for a certain distance. We
can thus assign, to any point $P\in%
\mathbb{R}
^{2}$, its signed distance $\delta(P)$ to $L_{0}$; $\delta(P)>0$ iff $P$ is in
the open half-plane bordered by $L_{0}$ into which $\nabla f$ points. A key
observation is that%
\begin{equation}
f(P)=\left\vert \nabla f\right\vert \cdot\delta(P)\text{.} \label{f(P)}%
\end{equation}
This is the principle \textquotedblleft rise$\,$=$\,$slope$\times
$run\textquotedblright\ -- the essence of linearity.

An entirely analogous situation takes place in the case of a three-variable
linear function $F(x,y,z)=ax+by+cz$, for which the level sets are the planes
normal to $\nabla F\equiv a\mathbf{i}+b\mathbf{j}+c\mathbf{k}$. For a point
$P\in%
\mathbb{R}
^{3}$, we similarly define its signed distance $\delta(P)$ to $F$'s level-$0 $
plane $\Pi_{0}$. Then,%
\begin{equation}
F(P)=\left\vert \nabla F\right\vert \cdot\delta(P)\text{.} \label{F(P)}%
\end{equation}

Equations (\ref{f(P)}) and (\ref{F(P)}) hold the key to the ensuing applications.

\section{Calculating distances without any formulae}

What is the distance between the two lines $L_{1}:=\{(x,y)\mid3x+4y=1\}$ and
$L_{31}:=\{(x,y)\mid3x+4y=31\}$? Even if you remember a relevant formula, you
might prefer this one-liner: $L_{1}$ and $L_{31}$ are the level-$1$ and
level-$31$ sets of $f(x,y):=3x+4y$ and therefore the distance (the
\textquotedblleft run\textquotedblright) between them equals%
\[
\frac{\text{the rise}}{\text{the slope}}=\frac{31-1}{\left\vert \nabla
f\right\vert }=\frac{30}{5}=6\text{ .}
\]
Similarly, the distance between the two planes $\Pi_{1}:=\{(x,y,z)\mid
2x+y+2z=1\}$ and $\Pi_{31}:=\{(x,y,z)\mid2x+y+2z=31\}$ equals%
\[
\frac{31-1}{\left\vert \nabla(2x+y+2z)\right\vert }=\frac{30}{3}=10\text{.}
\]

What is the distance between a point and a plane, for example the point
$P_{0}:=(1,2,1)$ and the plane $\Pi_{0}:=\{(x,y,z)\mid2x+3y+6z=0\}$? Note that
$\Pi_{0}$ is the level-$0$ plane of $F(x,y,z):=2x+3y+6z$ and the distance
sought therefore equals%
\[
\frac{\text{the change}}{\text{the rate}}=\frac{\left\vert F(P_{0})-F(\Pi
_{0})\right\vert }{\left\vert \nabla F\right\vert }=\frac{14-0}{7}=2\text{ .}
\]

Indeed, the\ relevant distance formulae can be derived by the same arguments.

\section{Deducing the Cauchy-Schwartz inequality by inspection}

For two $n$-tuples of real numbers $(a_{i})_{i=1}^{n}$ and $(x_{i})_{i=1}^{n}
$, the Cauchy-Schwartz inequality asserts that $\left\vert \sum_{i=1}^{n}%
a_{i}x_{i}\right\vert \leq\sqrt{\sum_{i=1}^{n}a_{i}^{2}}\sqrt{\sum_{i=1}%
^{n}x_{i}^{2}}$. For its proof, it suffices to treat the special case in which
$\sum_{i=1}^{n}x_{i}^{2}=1$; for, if $\sum_{i=1}^{n}x_{i}^{2}=s^{2}$, then
$\sum_{i=1}^{n}\left(  x_{i}/s\right)  ^{2}=1$.

We motivate our proof through an example. Let $S^{2}=\{(x,y,z)\mid x^{2}%
+y^{2}+z^{2}=1\}$, the unit $2$-sphere; let $F(x,y,z)=2x+y+2z$. Consider this
familiar problem:\ Find $\max_{S^{2}}F$ and $\min_{S^{2}}F$. In light of
Eq.\thinspace(\ref{F(P)}), the solution can be argued by mere inspection,
almost without words, as follows.

View the level-$0$ plane $\Pi_{0}$ of $F$ as $S^{2}$'s equatorial plane and
its normal direction $2\mathbf{i}+\mathbf{j}+2\mathbf{k}$ as North. By
Eq.\thinspace(\ref{F(P)}), the extrema of $F$ on $S^{2}$ are attained
precisely at the points on $S^{2}$ that are farthest from $\Pi_{0}$, i.e., the
North Pole $P_{+}$ and the South Pole $P_{-}$. Clearly, $P_{\pm}=\left(
\pm\frac{2}{3},\pm\frac{1}{3},\pm\frac{2}{3}\right)  $, at which
Eq.\thinspace(\ref{F(P)}) again gives the $F$-value: $F(P_{\pm})=\left\vert
\nabla F\right\vert \cdot\delta(P_{\pm})=3\cdot(\pm1)$. Hence, $\left\vert
2x+y+2z\right\vert \leq3$ for $(x,y,z)\in S^{2}$.

This argument applied to the linear function $ax+by+cz$ yields that%
\[
\left\vert ax+by+cz\right\vert \leq\left\vert \nabla(ax+by+cz)\right\vert
=\sqrt{a^{2}+b^{2}+c^{2}}\text{ \quad for }(x,y,z)\in S^{2}\text{ ,}%
\]
which, as remarked earlier, implies the Cauchy-Schwartz inequality for
$n=3$.\begin{figure}[h]
\centering\includegraphics[scale=0.6]{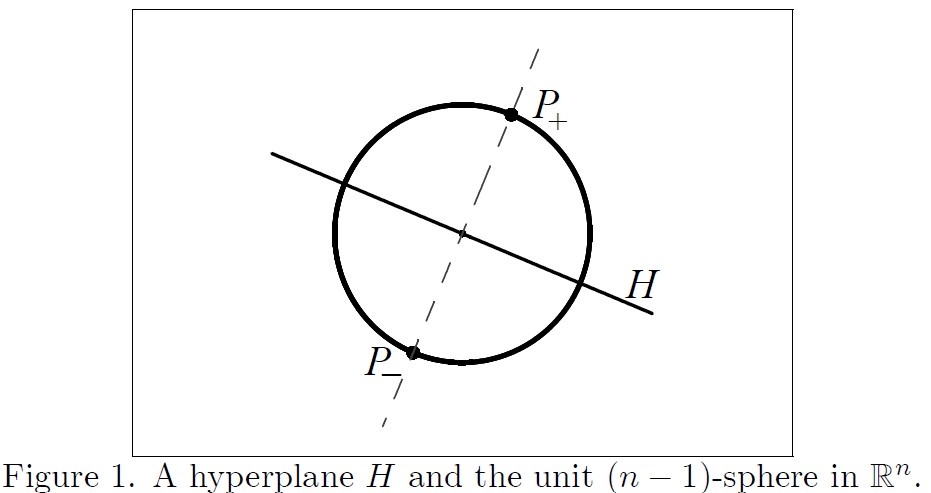}\end{figure}

The preceding argument even works in higher dimensions. Let $G$ be a
homogeneous degree-1 polynomial in $x_{1},\cdots,x_{n}$. Its level-$0$ set $H$
is a hyperplane (a 1-codimensional subspace) in $%
\mathbb{R}
^{n}$. By an orthogonal transformation of coordinates, we may install a new
coordinate system $(x_{1}^{\prime},\cdots,x_{n}^{\prime})$ such that $H$ is
the subspace $x_{n}^{\prime}=0$ and that, for $P\in%
\mathbb{R}
^{n}$, its primed $n$th coordinate $x_{n}^{\prime}(P)$ is its signed distance
to $H$. See Figure 1. We thus have%
\[
G(P)=\left\vert \nabla G\right\vert \cdot x_{n}^{\prime}(P).
\]
The unit $(n-1)$-sphere $S^{n-1}$ in $%
\mathbb{R}
^{n}$ satisfies the same equation in the primed coordinates: $\sum_{i=1}%
^{n}\left(  x_{i}^{\prime}\right)  ^{2}=1$, from which it is clear that
$x_{n}^{\prime}$ on $S^{n-1}$ attains its extrema at precisely the two points
$P_{\pm}$ whose primed coordinates are $(0,\cdots,0,\pm1)$, implying that
$\left\vert G(P)\right\vert \leq\left\vert \nabla G\right\vert $ for $P\in
S^{n-1}$.

\section{Visualizing the GM-AM inequality}

For $x_{1},\cdots,x_{n}>0$, the GM-AM inequality asserts that their geometric
mean $\sqrt[n]{x_{1}\cdots x_{n}}$ is no greater than their arithmetic mean
$\left(  x_{1}+\cdots+x_{n}\right)  /n$. The following special case is an
equivalent statement:

\begin{center}
$x_{1}+\cdots+x_{n}\geq n$\quad for $x_{1},\cdots,x_{n}>0$ such that
$\prod\nolimits_{i=1}^{n}x_{i}=1$ .
\end{center}

\noindent For, if $\prod\nolimits_{i=1}^{n}x_{i}=p$, then $\prod
\nolimits_{i=1}^{n}\left(  x_{i}/\sqrt[n]{p}\right)  =1$. We show a visual
argument for this statement in the case $n=2$ and suggest one for
$n=3$.\begin{figure}[h]
\centering\includegraphics[scale=0.4]{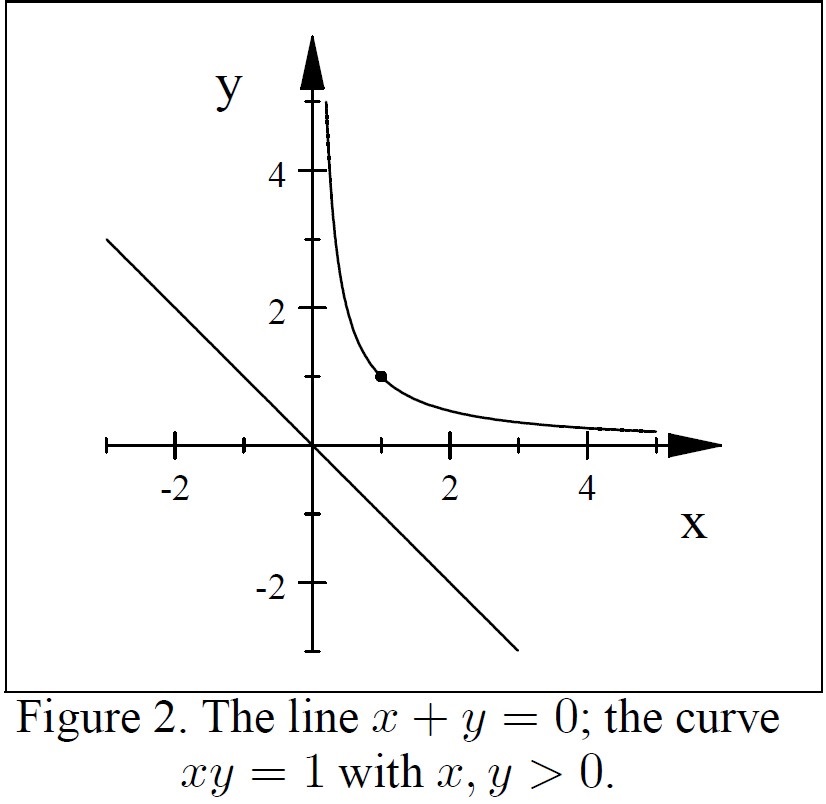}\end{figure}\begin{figure}[h]
\centering\includegraphics[scale=0.4]{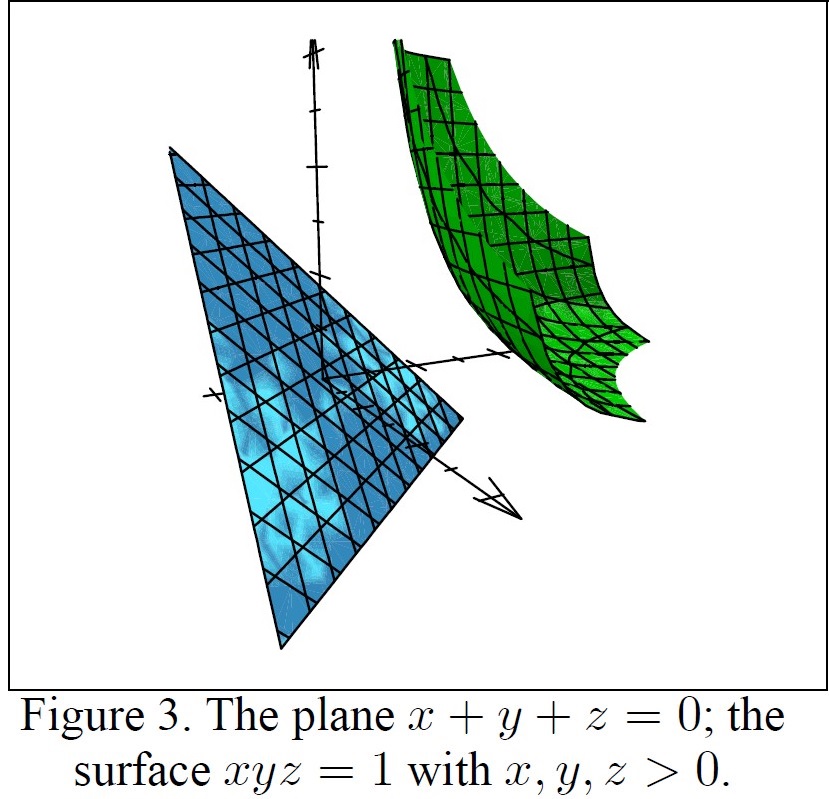}\end{figure}

Let $C=\{(x,y)\mid x,y>0;\;xy=1\}$. Where on $C$ does $x+y$ attain its least
value? The answer is given by Eq.\thinspace(\ref{f(P)}):\ at those points on
$C$ that are closest to the line $L_{0}=\{(x,y)\mid x+y=0\}$. It is apparent
from Figure 2 that $(1,1)$ is the only minimizer of this distance and
therefore $x+y\geq1+1$ for $(x,y)\in C$. For proof, view $C$ as the graph of a
\textit{convex} function $h$ defined on $L_{0}$. (In a new coordinate system
making a $45^{\circ}$-angle with the standard one, an equation for $C$ in the
new coordinates $(x^{\prime},y^{\prime})$ is $y^{\prime}=\sqrt{1+(x^{\prime
})^{2}}$.) As $h$ is convex, it attains its global minimum exactly once; see
\cite{Flemming} for convex functions. Symmetry then dictates where the
minimizer is.

For the case $n=3$, let $S=\{(x,y,z)\mid x,y,z>0;\;xyz=1\}$, which is the
graph of the function $g(x,y)=(xy)^{-1}$ on the first quadrant $(0,\infty
)^{2}\subset%
\mathbb{R}
^{2}$. On $S$, the quantity $x+y+z$ attains its least value at a point closest
to the plane $\Pi_{0}=\{(x,y,z)\mid x+y+z=0\}$ (by (\ref{F(P)})). Looking at
Figure 3, it is plausible that $(1,1,1)$ is the uniqe minimizer of this
distance, which then implies that $x+y+z\geq1+1+1$ for $(x,y,z)\in S$. In
fact, this situation is almost identical to the previous case: $S$ can be
viewed as \textit{the graph of a convex function on }$\Pi_{0}$; convexity then
entails the existence and uniqueness of a minimizer and symmetry dictates its
location. To make this visual argument rigorous, we must analytically check
that $S$ can indeed be viewed as the graph of a function on $\Pi_{0}$ and that
this function is indeed convex. For the first claim, we need to show that $S$
passes the \textquotedblleft vertical line test\textquotedblright\ where a
\textquotedblleft vertical line\textquotedblright\ means a line perpendicular
to the domain $\Pi_{0}$. From a point $(x,y,-x-y)\in\Pi_{0}$, proceed in the
\textquotedblleft vertical\textquotedblright\ direction $\nabla
(x+y+z)=\mathbf{e}_{1}+\mathbf{e}_{2}+\mathbf{e}_{3}$ to arrive at points of
the form%
\[
P(t)=(x+t,\,y+t,\,-x-y+t)\text{.}%
\]
For how many values of $t$ will $P(t)\in S$?\ A computer can check for us that
the cubic equation in $t$%
\[
(x+t)(\,y+t)(-x-y+t)=1
\]
has exactly one real solution. For the other claim concerning convexity, first
note that the function $g(x,y)=(xy)^{-1}$ is convex on the first quadrant; one
can verify this by checking that the Hessian quadratic form of $g$ is positive
definite at every point in the first quadrant. Hence $S$ lies on one side of
each of its tangent planes, a fact that is invariant under a change of perspective.

\end{document}